\documentclass[a4paper]{article}       

\usepackage{amsmath,mathtools,amsthm}
\usepackage{amssymb, authblk}
\usepackage[mathscr]{eucal}
\usepackage{bm}
\usepackage{bbm}
\usepackage[shortlabels]{enumitem}
\usepackage{hyperref}
\usepackage[usenames,dvipsnames]{xcolor}

\usepackage{tikz}
\usepackage{graphicx, caption,subcaption}
\usepackage{todonotes}

\newcommand{\prob}{\mathbb{P}}
\newcommand{\Prob}[1]{\prob\left(#1\right)}

\newcommand{\expec}{\mathbb{E}}
\newcommand{\Exp}[1]{\expec\left[#1\right]}

\allowdisplaybreaks

\newtheorem{theorem}{Theorem}[section]

\newtheorem{proposition}[theorem]{Proposition}

\numberwithin{equation}{section}

\usepackage{a4wide}
\title{Infinite-horizon Fuk-Nagaev inequalities}
\author{A.J.E.M. Janssen$^1$ and B. Zwart$^{1,2}$}
\date{%
  {\small  $^1$Eindhoven University of Technology\\%
    $^2$CWI Amsterdam} \\[2ex]%
    \today
}

\begin{document}
\maketitle
 \begin{abstract}
     We develop explicit bounds for the tail of the distribution for the all-time supremum of a random walk with negative drift, where the increments have a truncated heavy-tailed distribution. As an application, we consider a ruin problem in the presence of re-insurance.

\vspace{1mm}

 \end{abstract}
\section{Introduction and main results}

Let $X_i, i\geq 1,$ be a real-valued i.i.d.\ sequence with $\mu = \Exp{X_1}<0$, $\mu_\beta=\Exp{|X_1|^\beta}<\infty$ for some $\beta>1$, and $\Exp{e^{sX_1}}=\infty$ for $s>0$, so that $X_1$ is heavy-tailed. Define for $y>0$, $S_0(y)=0, S_n(y) = \sum_{i=1}^n \min \{X_i, y\}, n\geq 1,$ and
\begin{equation}
\label{protagonist}
    M(y) = \sup_{n\geq 0} S_n(y).
\end{equation}
We derive explicit bounds for $\Prob{M(y)>x}$, motivated by our desire to understand the impact of truncation on the behavior of $M(y)$. Truncated heavy tails naturally occur in risk theory \cite{AsmussenPihlsgaard}, where large claims can be re-insured. Other applications can be found in queueing theory, in which large jobs may get terminated \cite{Jelenkovic99}, and scale-free random graphs, where the degree of a vertex cannot be larger than the size of the graph \cite{kerriou2023fewestbigjumps,  Stegehuis23}.

Asymptotic properties of $\Prob{S_n(y)> x}$ where $x$ and $y$ go to $\infty$ at a rate comparable to $n$ have been analyzed by Chakrabarty \cite{Chakra12, Chakra17}.
Asymptotics for the finite-time ruin probability for various re-insurance contracts have been derived in \cite{albrecher2020finite, bazhba2020sample,
chen2019efficient, rhee2016arxiv, dombry2022hidden}.

For $M(y)$, such results seem unavailable.
If $y$ is fixed, Cramér-Lundberg theory \cite{asmussen2010ruin} implies that
\begin{equation}
\label{eq:cl}
    \Prob{M(y)> x} \leq e^{-\gamma(y)x},
\end{equation}
with $\gamma(y)$ the unique strictly positive solution of the Cramér-Lundberg equation
\begin{equation}
\label{fpe}
   \Exp{e^{s \min \{X_1,y\}}}=1,
\end{equation}
so that $\gamma(y) = \sup \{s\geq 0:\Exp{e^{s \min \{X_1,y\}}}\leq 1\}$.
Since $X_1$ is heavy-tailed, $\gamma(y)\rightarrow 0$ as $y\rightarrow\infty$. In
\cite{AsmussenPihlsgaard}, the authors derive asymptotic expansions for $\gamma(y)$ in the slightly different setting of a L\'evy process with truncated jumps;
their results correspond to the case where $\Prob{X_1>x}$ is either regularly varying or of the form $e^{-\sqrt{x}}$.

We mainly focus on non-asymptotic estimates. In particular, we aim to derive upper bounds for $\Prob{M(y)>x}$. Such bounds can be valuable, as convergence rates of asymptotic estimates can be slow for heavy-tailed random variables \cite{mikosch1998large}.
Our first main result is the following theorem.


\begin{theorem}
\label{thm-main}
    There exists a $y_\beta<\infty$ such that
    \begin{equation}
    \label{eq-main-beta}
        \Prob{M(y)>x} \leq y^{-(\beta-1) x/y}, \hspace{1cm}  x>0, y>y_\beta.
    \end{equation}
\end{theorem}
An explicit expression for $y_\beta$ is given in (\ref{def-ybeta}) below.
Theorem \ref{thm-main} is aimed at random walks with step size distributions for which $\Exp{X_1^k}=\infty$ for some $k>0$.
Our next result covers cases where $\Prob{X_1>x}$ has a lighter tail, e.g., $\Prob{X_1>x}$ is of the form $\exp \{ -(\log x)^\xi\}, \xi >1$ (which includes log-normal tails) or
$\exp \{ -x^\xi\}, \xi \in (0,1)$ (which includes Weibull tails). Define
\[
q(x) = -\log \Prob{X_1>x}.
\]
\begin{theorem}
\label{thm-main2}
Assume that $\mu_2<\infty$ and suppose that there exists a $\kappa \in (0,1)$ and $y_\kappa>0$  such that
$(\log y)^{1+\kappa} \leq q(y) \leq y^{1-\kappa}$ for $y\geq y_\kappa$, and that $q(y)$ is concave.
 Then, for every $\eta<1$, there exists a $y_\eta^*<\infty$ such that
  \begin{equation}
  \label{eq-logconvex}
        \Prob{M(y)>x} \leq \left( \frac {y\Prob{X_1>y}}{\eta |\mu|} \right)^{\frac xy}, \hspace{1cm}  x>0, y>y^*_\eta.
    \end{equation}
\end{theorem}
An explicit expression for $y^*_\eta$  is given in (\ref{def-yeta}) below. 

Theorem \ref{thm-main} and Theorem \ref{thm-main2} complement similar bounds for $\Prob{S_n(y)>x}$, which are known as Fuk-Nagaev inequalities after \cite{FN71}.
The approach to proving such inequalities, surveyed in \cite{Nagaev1979}, has been modified in \cite{Borovkov1972} and \cite{Fuk1973} to investigate $\Prob{ \max_{k\leq n} S_k(y) > x}$; our results do not seem to follow from the estimates in \cite{Borovkov1972, Fuk1973}.  Related bounds for $M(\infty)$ have been derived in \cite{Kugler_Wachtel_2013}, which are shown to be sharp (for large $x$ and/or $\mu$ small) for distributions with power-law tails.

A key step in the approach in \cite{Nagaev1979} and the other cited works includes using the Chernoff bound for $\Prob{S_n(y)>x}$ and deriving
convenient upper bounds for $\Exp{e^{s\min \{X, y\}}}$.
Our proof ideas are similar, but we rely on (\ref{eq:cl}) and (\ref{fpe}) instead, and derive lower bounds for $\gamma(y)$.
We have tried to strike a balance between the accuracy of such lower bounds, and user-friendliness (in terms of explicitness)
of the final expressions (\ref{eq-main-beta}), (\ref{eq-logconvex}). This is partly made possible by restricting the range of $y$ over which these bounds apply.
It is possible to extend this range, at the expense of making the bound less sharp for larger values of $y$, using the explicit form of $y_\beta$ and $y_\eta^*$
given in (\ref{def-ybeta}) and (\ref{def-yeta}) below. For more details, we refer to Sections \ref{sec-proofmaintheorem} and \ref{sec-proofmaintheorem2}.
For a general framework that can be applied to develop sharp estimates for $\Exp{e^{s\min \{X, y\}}}$, we refer to \cite{Janssen2024}.

To illustrate the applicability and accuracy of our bounds, we consider a variation of the Sparre-Andersen risk model where claims to exceed a value
$y$ are re-insured.
We analyze the infinite-time ruin probability $P_a(x)$ if the initial capital is $x$, and $y$ is chosen as $y=ax$, where $a$ is a fixed constant and $x\rightarrow\infty$, using Theorem \ref{thm-main}.
Finite-time ruin versions of this re-insurance problem have been analyzed in \cite{bazhba2020sample, rhee2016arxiv} using
sample-path large deviation principles for heavy tails.
Extending such a result to the infinite horizon case requires an additional interchange of limit argument. In Section \ref{sec-proofcor} below, we illustrate how our bounds can help in this regard.

The result in Section \ref{sec-proofcor} shows that if $y=ax$, our bounds almost predict the correct asymptotic behavior,
except for certain round-offs. Informally, if claims sizes have a tail of the form $x^{-\alpha}$, ruin occurs due to the consequence of $\lceil 1/a \rceil$ large claims, leading to a ruin probability which behaves like $x^{-(\alpha-1)\lceil 1/a \rceil}$, while Theorem \ref{thm-main} would suggest
a behavior of $x^{-(\alpha-1)/a}$.

The remainder of this article is organized as follows. In Section \ref{sec-proofmaintheorem}, we give a proof of Theorem \ref{thm-main}.
The proof of Theorem \ref{thm-main2} is given in Section \ref{sec-proofmaintheorem2}. The application to a risk model is given
in Section \ref{sec-proofcor}.
The appendix presents the results from \cite{Janssen2024} that are needed in the main body of the paper.

\newpage

\section{Proof of Theorem \ref{thm-main}}
\label{sec-proofmaintheorem}

We use the following notation: $x^+ = \max \{x,0\}$, $F(x) = \Prob{X_1>x}$, $\bar F(x)=1-F(x)$. Indicator functions are written as $I(\cdot)$.
Denote
$\mu_\beta^+ = \Exp{(X_1^+)^\beta}$,
$\mu_\beta^- = \Exp{((-X_1)^+)^\beta}$.
We drop $\beta$ from these expressions if $\beta=1$, for example, $\mu^- = \Exp{ (-X_1)^+ }$.
We write, for $\beta>1$, $\mu_\beta= \Exp{ |X_1|^\beta}$. Observe that $\mu_\beta =\mu_\beta^++\mu_\beta^-$.
When we apply this identity, we always take $\beta>1$ and warn that, in contrast, $\mu = \mu^+-\mu^-$.

%
%

In view of the inequality (\ref{eq:cl}), Equation (\ref{eq-main-beta}) holds if we can show that there exists a $y_\beta$ such that $\gamma(y) \geq s(y)=((\beta-1) \log y)/y$ for $y>y_\beta$.
For this, using the equation (\ref{fpe}) for $\gamma(y)$, it is sufficient to show that $ \Exp{e^{s(y) \min \{X_1,y\}}}\leq 1$.
We therefore derive an appropriate bound for $\Exp{e^{s \min \{X_1, y\}}}$.
Let $\delta = \min \{2,\beta\}$. We have, for $s>0$,
\begin{equation}
    \label{2.2}
     \Exp{e^{s \min \{X_1,y\}}} = \int_{-\infty}^0e^{sx} dF(x) +\int_{0}^ye^{sx} dG(x),
\end{equation}
with $G(x)=\Prob{\min \{X_1, y\} \leq x}$.
We estimate
\begin{align}
    \int_{-\infty}^0 e^{sx} dF(x) &=  \int_{-\infty}^0  (e^{sx}-1-sx) dF(x) + \int_{-\infty}^0 (1+sx) dF(x) \nonumber \\
    &= \int_{-\infty}^0  \frac{e^{sx}-1-sx}{(-sx)^\delta} s^\delta (-x)^\delta dF(x) + F(0) - s\mu^- \nonumber \\
    &\leq s^\delta MG_\delta \mu_\delta^- + F(0)-s\mu^-, \label{2.4}
\end{align}
where (see the Appendix)
\begin{equation}
\label{2.5}
    MG_\delta = \max_{u\geq 0} \frac{e^{-u}-1+u}{u^\delta} \leq \left(\frac{\delta-1}{\delta}\right)^{\delta-1}(2-\delta)^{2-\delta} \leq \frac 1\delta.
\end{equation}
Next, we bound the integral $\int_{0}^ye^{sx} dG(x)$ using arguments similar to those in the proof of Lemma 1.4 of \cite{Nagaev1979}. Observe that
\begin{align}
    \int_{0}^ye^{sx} dG(x) &= \int_{0}^y(e^{sx}-1-sx) dG(x) + \int_{0}^y(1+sx) dG(x)\nonumber \\
    &\leq \int_{0}^y(e^{sx}-1-sx) dG(x) + \bar F(0) + s\mu^+.  \label{2.9}
\end{align}
To bound the remaining integral at the right-hand side of (\ref{2.9}), we distinguish between the cases that $0\leq s\leq \beta/y$ and $s> \beta/y$. In the former case, we have
\begin{equation}
    \label{2.10}
    \int_{0}^y(e^{sx}-1-sx) dG(x)  \leq \int_{0}^{\beta/s}(e^{sx}-1-sx) dG(x).
\end{equation}
Recalling $\delta = \min \{2,\beta\}$, we have for $0\leq x\leq \beta/s$,
\begin{equation}
\label{2.11}
    0\leq \frac{e^{sx}-1-sx}{(sx)^\delta} = e^{sx} \frac{1-(1+sx)e^{-sx}}{(sx)^\delta}
    \leq e^{sx} ME_\delta \leq e^\beta ME_\delta,
\end{equation}
where, see the Appendix,
\begin{equation}
\label{2.12}
    ME_\delta = \max_{v\geq 0} \frac{1-(1+v)e^{-v}}{v^\delta} \leq    \frac{(2-\delta)^{2-\delta}}{4-\delta} \leq \frac 1{1+\frac 12 \delta}.
\end{equation}
Hence, for the case $0\leq s\leq \beta/y$, we have
\begin{equation}
    \label{2.13}
    \int_{0}^y(e^{sx}-1-sx) dG(x)  \leq e^\beta ME_\delta \int_0^{\beta/s} (sx)^\delta d G(x) \leq e^\beta s^\delta ME_\delta \mu_\delta^+.
\end{equation}
If $s> \beta/y$, we have from (\ref{2.13}),
\begin{equation}
\label{2.14}
    \int_{0}^y(e^{sx}-1-sx) dG(x) \leq e^\beta s^\delta ME_\delta \mu_\delta^+ + \int_{\beta/s}^y (e^{sx}-1-sx) dG(x).
\end{equation}
Now, observe that
\begin{equation}
\label{2.15}
    \frac{d}{dx} \left[ \frac{e^{sx}-1-sx}{x^\beta}\right] = x^{-\beta-1} (e^{sx}(sx-\beta)+(\beta-1)sx + \beta)>0
\end{equation}
for $x\geq \beta/s$, and so
\begin{equation}
    \label{2.16}
\int_{\beta/s}^y (e^{sx}-1-sx) dG(x) = \int_{\beta/s}^y \frac{e^{sx}-1-sx}{x^\beta} x^\beta dG(x) \leq  \frac{e^{sy}-1-sy}{y^\beta} \mu_\beta^+.
\end{equation}
We conclude that, when $s\geq \beta/y$,
\begin{equation}
    \label{2.17}
\int_{0}^y(e^{sx}-1-sx) dG(x) \leq e^\beta s^\delta ME_\delta \mu_\delta^+ + \frac{e^{sy}-1-sy}{y^\beta} \mu_\beta^+.
\end{equation}
Evidently (see (\ref{2.13})), the bound (\ref{2.17}) also holds when $0\leq s\leq \beta/y$.
Returning to (\ref{2.2}) and using (\ref{2.4}), (\ref{2.9}), (\ref{2.13}), (\ref{2.17}), we get
\begin{equation}
\Exp{e^{s \min \{X_1, y\}}} \leq s^\delta MG_\delta \mu_\delta^-+F(0)-s\mu^-+\bar F(0) + s\mu^++  e^\beta s^\delta ME_\delta \mu_\delta^+ + \frac{e^{sy}-1-sy}{y^\beta} \mu_\beta^+.  \label{2.18}
\end{equation}
Using that $F(0)+\bar F(0)=1$ and $\mu^+-\mu^-=\mu$, we see that
\begin{align}
    \Exp{e^{s \min \{X_1, y\}}} &\leq
1+s\mu+ s^\delta MG_\delta \mu_\delta^-+  e^\beta s^\delta ME_\delta \mu_\delta^+ + \frac{e^{sy}-1-sy}{y^\beta} \mu_\beta^+ \nonumber\\
&= 1+s\left(\mu+s^{\delta-1} MG_\delta \mu_\delta^-+  e^\beta s^{\delta-1} ME_\delta \mu_\delta^+ + \frac{e^{sy}-1-sy}{sy^\beta} \mu_\beta^+\right). \label{2.19}
\end{align}
We consider this with $s=s(y)= \frac 1y (\beta-1) \log y$, so that $e^{sy}=y^{\beta-1}$.
For $y>1$ we therefore obtain,
\begin{equation}
    \label{2.20}
    0\leq \frac{e^{sy}-1-sy}{sy^\beta} = \frac{y^{\beta-1}-1-sy}{y^{\beta-1}(\beta-1)\log y} \leq \frac 1{(\beta-1) \log y}.
\end{equation}
Concluding, we get for $y>1$ and $s= \frac 1y (\beta-1) \log y$,
\begin{equation}
      \label{2.21}
      \Exp{e^{s \min \{X_1, y\}}} \leq
      1+s\left(\mu+s^{\delta-1} (MG_\delta \mu_\delta^- + e^\beta  ME_\delta \mu_\delta^+) + \frac{\mu_\beta^+}{(\beta-1)\log y}\right).
\end{equation}
The quantity on the right-hand side between brackets is negative for large enough $y$, since $\mu<0$. To quantify this statement, write
\begin{equation}
      \label{2.22}
  K = \mu_\delta^- MG_\delta  +    e^\beta  ME_\delta \mu_\delta^+, \hspace{0.5cm} L= \mu_\beta^+,
\end{equation}
so that the relevant quantity in (\ref{2.21}) takes the form
\begin{equation}
\label{2.23}
    Q = \mu + \left(Ks^{\delta-1}+ \frac{L}{(\beta-1) \log y} \right) = \mu + \frac{L+(\beta-1)K s^{\delta-1} \log y}{(\beta-1) \log y}.
\end{equation}
With $y>1$ and $s= \frac 1y (\beta-1) \log y$, observe that
\begin{equation}
    \label{2.24}
    (\beta-1) s^{\delta-1} \log y = (\beta-1)^\delta \frac{(\log y)^\delta}{y^{\delta-1}} \leq C (\beta-1)^\delta,
\end{equation}
where
\begin{equation}
\label{2.25}
    C=\max_{y>1} \frac{(\log y)^\delta}{y^{\delta-1}} = \left(\max_{y>1} \frac{\log y}{y^{\frac{\delta-1}{\delta}}}\right)^\delta = \left(\frac{\delta}{(\delta-1)e}\right)^\delta,
\end{equation}
with maximizing $y$ given by $\exp\{\delta/(\delta-1)\}$. Thus,
\begin{equation}
\label{2.26}
    (\beta-1) s^{\delta-1} \log y \leq \left( \frac{\beta-1}{\delta-1}\frac \delta e\right)^\delta.
\end{equation}
Consequently, the quantity $Q$ in (\ref{2.23}) is negative when
\begin{equation}
    Q \leq \mu + \frac{L+K\left( \frac{\beta-1}{\delta-1}\frac \delta e\right)^\delta}{(\beta-1) \log y}.
\end{equation}
Hence, the quantity $Q$ in (\ref{2.23}) is negative when
\begin{equation}
\label{def-ybeta}
    y>y_\beta = \exp\left\{\frac {L+K\left(\frac{\beta-1}{\delta-1}\frac \delta e\right)^\delta}{(\beta-1) |\mu|}\right\}.
\end{equation}
Here, $K$ and $L$ are given in (\ref{2.22}) where $MG_\delta$ and $ME_\delta$ can be bounded using (\ref{2.5}) and (\ref{2.12}).

Thus, for $y>y_\beta$, it follows that $\gamma(y) \geq \frac 1y (\beta-1) \log y$.
Inserting this bound into (\ref{eq:cl}) completes the proof
 of (\ref{eq-main-beta}).\\

\noindent
{\bf Remarks.}
\begin{enumerate}
    \item Upon inspection of the proof, at the expense of increasing the threshold $y_\beta$, we could have chosen
$s=s(y) = ((\beta-1) \log y + \log\log y + \xi)/y$ with $e^\xi < (\beta-1) |\mu| / \mu_\beta^+$, which would lead to a slightly sharper bound
of the form
\begin{equation}
    \label{eq-main-beta2}
        \Prob{M(y)>x} \leq y^{-(\beta-1)\frac xy} (\log y)^{-\frac xy} e^{-\gamma \frac xy}   , \hspace{1cm}  x>0, y>y_\beta.
    \end{equation}
\item If one wishes an upper bound valid for all $y$, one may simply use
\begin{equation}
    \Prob{M(y)>x} \leq (y/y_\beta)^{-(\beta-1)x/y}.
\end{equation}
\end{enumerate}

\section{Proof of Theorem \ref{thm-main2}}
\label{sec-proofmaintheorem2}

As in the proof of Theorem \ref{thm-main}, we bound $\Exp{e^{s \min \{X_1,y\}}}$,
now using the assumptions made in the formulation of Theorem \ref{thm-main2}.
For a fixed $\eta \in (0,1)$, we then make a particular choice for $s=s(y)$, and we construct a $y_\eta^*$ such that $\Exp{e^{s \min \{X_1,y\}}}\leq 1$ for $y\geq y_\eta^*$ and $s=s(y)$. By Cram\'er-Lundberg theory, this implies $\gamma(y)\geq s(y)$, and inserting this inequality into (\ref{eq:cl}), we obtain (\ref{eq-logconvex}).

Let $s>0$ and $y>1/s$. We have
\begin{align}
    \Exp{e^{s \min \{X_1,y\}}} &= \int_{-\infty}^ye^{sx} dF(x) + e^{sy} \bar F(y) \nonumber \\
    &= \int_{-\infty}^0e^{sx} dF(x) + \bar F(0) +\int_0^{1/s} se^{sx} \bar F(x) dx +\int_{1/s}^y se^{sx} \bar F(x) dx. \label{3.18}
 \end{align}
We bound the three integrals in (\ref{3.18}). For the first integral, we can use (\ref{2.4}), where we observe that $\delta = \min \{2,\beta \}=2$ since $\mu_2<\infty$,
and that $MG_2=1/2$. Thus, we have
\begin{equation}
    \label{3.19}
    \int_{-\infty}^0e^{sx} dF(x) \leq \frac{1}{2}s^2 \mu_2^-+F(0)-s\mu^-.
\end{equation}
To bound the second integral in (\ref{3.18}), we use the inequality $e^z\leq 1+ze^z$ with $z=sx$ and the inequalities
\begin{equation}
\label{3.20}
    \int_0^u \bar F(x)dx \leq \mu^+, \hspace{0.5cm} \int_0^u 2x\bar F(x)dx \leq \mu_2^+, \hspace{0.5cm} u\geq 0,
\end{equation}
to obtain
\begin{align}
    \int_0^{1/s} se^{sx} \bar F(x) dx &\leq \int_0^{1/s} s(1+sxe^{sx}) \bar F(x) dx
    \leq s \int_0^{1/s} \bar F(x) dx + \frac{1}{2} s^2 \int_0^{1/s} 2x e \bar F(x) dx \nonumber \\
    &\leq s\mu^+ +s^2 \frac e2 \mu_2^+. \label{3.21}
\end{align}
To bound the third integral in (\ref{3.18}),
note that $sx-q(x)$ is convex, so that
\begin{equation}
    sx-q(x) \leq \max \{1-q(1/s), sy-q(y)\}, \hspace{1cm} x\in [1/s,y].
    \label{3.22}
\end{equation}
Consequently,
\begin{equation}
    \int_{1/s}^y se^{sx} \bar F(x) dx=\int_{1/s}^y se^{sx-q(x)} dx \leq sy \max \{e^{1-q(1/s)}, e ^{sy-q(y)}\}. \label{3.23}
\end{equation}
From (\ref{3.19}), (\ref{3.21}) and (\ref{3.23}), we then see that
\begin{align}
     \Exp{e^{s \min \{X_1,y\}}} &\leq (\frac{1}{2}s^2 \mu_2^-+F(0)-s\mu^-) + \bar F(0)+
     (s\mu^+ +s^2 \frac e2 \mu_2^+)+ sy \max \{e^{1-q(1/s)}, e ^{sy-q(y)}\} \nonumber \\
     &= 1+ s\left(\mu+  s \frac e2 \mu_2^+ + \frac 12 s\mu_2^- + y \max \{e^{1-q(1/s)}, e ^{sy-q(y)}\} \right) \nonumber \\
     &\leq  1+ s\left(\mu+  s \frac e2 \mu_2 + y \max \{e^{1-q(1/s)}, e ^{sy-q(y)}\} \right).\label{3.24}
\end{align}
In the last display, we used that $\mu^+-\mu^-=\mu, F(0)+\bar F(0)=1$, and $\mu_2^++\mu_2^-=\mu_2$.
We now want to choose $s(y)$ such that the expression
\begin{equation}
\label{3.26}
    Q= \mu+  s \frac e2 \mu_2 + y \max \{e^{1-q(1/s)}, e^{sy-q(y)}\} \leq 0,
\end{equation}
when $y$ is sufficiently large.
We set
\begin{equation}
    \label{3.27}
    s=s(y):=\frac{1}{y}(q(y)-\log y + r), \hspace{0.5cm} r= \log (|\mu| \eta),
\end{equation}
with $\eta\in (0,1)$ fixed. Then, we have
\begin{equation}
   \label{3.28}
    Q= \mu+  s \frac e2 \mu_2 +  \max \{y e^{1-q(1/s)}, e^{r}\}.
\end{equation}
We now invoke the assumption made in Theorem \ref{thm-main2} that there exists a $\kappa \in (0,1)$ and $y_\kappa>0$ such that
\begin{equation}
\label{3.29}
    (\log y)^{1+\kappa} \leq q(y) \leq y^{1-\kappa}, \hspace{0.5cm} y\geq y_\kappa.
\end{equation}
We have
\begin{equation}
\label{3.30}
    y \exp \{- q(1/s(y))\} \leq y \exp \{-q(y^\kappa)\},
\end{equation}
when $y\geq \max \{ e^r, y_\kappa\}$. Indeed, by monotonicity of $q$, the inequality in (\ref{3.30}) is equivalent to
\begin{equation}
    \frac 1{s(y)} = \frac y{q(y)-\log y + r} \geq y^\kappa,
\end{equation}
and this follows from the second inequality in (\ref{3.29}) and $y\geq e^r$. Furthermore, from the first inequality in (\ref{3.29}), we have
\begin{equation}
    \label{3.32}
    q(y^\kappa) \geq (\log (y^\kappa))^{1+\kappa}=(\kappa \log y)^{1+\kappa}, \hspace{0.5cm} y\geq y_\kappa^{1/\kappa}.
\end{equation}
Now, define
\begin{equation}
    \label{3.33}
    y_r = \max \left\{ \sup \{y: ye^{-(\kappa \ln y)^{\kappa+1}} \geq e^{r-1}\}, e^r, y_\kappa^{1/\kappa} \right\}.
\end{equation}
Then, by (\ref{3.30}), (\ref{3.32}), and the definition of $y_r$ in (3.33), for $y\geq y_r$,
\begin{equation}
    \max \{y e^{1-q(1/s)}, e^{r}\} \leq \max \{y e^{1- (\kappa \ln y)^{\kappa+1} }, e^{r}\} \leq e^r.
\end{equation}
Therefore, $Q$ in (\ref{3.26}) and (\ref{3.28}) is bounded by
\begin{equation}
    \label{3.35}
    Q \leq \mu + \frac 12 es\mu_2 + e^r \leq \mu + \frac 12 e\mu_2 y^{-\kappa}+ e^r, \hspace{0.5cm} y\geq y_r,
\end{equation}
where the latter inequality follows from the definition of $s=s(y)$ in (\ref{3.27}) and the second inequality in (\ref{3.29}).
Since $e^r = \eta |\mu|$, we conclude that $Q\leq 0$ when $y\geq y_\eta^*$, where
\begin{equation}
\label{def-yeta}
    y_\eta^* = \max \left\{ y_r, \left(\frac{e\mu_2/2}{|\mu|(1-\eta)}\right)^{1/\kappa} \right\}.
\end{equation}
We conclude that $s(y) \leq \gamma (y)$ for $y>y_\eta^*$. The proof of Theorem \ref{thm-main2} now follows by combining the inequalities  $s(y) \leq \gamma (y)$ and (\ref{eq:cl}), and using the definition of $s(y)$ in (\ref{3.27}).\\

\noindent
{\bf Remarks.}
\begin{enumerate}
    \item The value of $y_\eta^*$ has not been optimized. In particular cases, the structure of $q(y)$ may be exploited to get sharper estimates, especially if $q(y)$ is of the form $(\log y)^{1+\kappa}$ or $y^{1-\kappa}$.
    \item Our assumptions on $q(y)$ may be slightly generalized. However, we conjecture that (\ref{eq-logconvex}) does not hold when the distribution of $X_1$ is too close to the exponential distribution, in particular if $q(y)= y/(\log y)^\zeta$, $\zeta\in (0,1)$.
\end{enumerate}

\section{Application to a re-insurance risk model}
\label{sec-proofcor}

Consider an insurance problem where premiums come in at rate $c$, claim sizes $B_i,i\geq 1$ are i.i.d.\, and the time between claims is governed by a renewal process, with $A_i, i\geq 1$, inter-renewal times which are independent of the claim sizes. Claim sizes exceeding a value of $y=ax$ are covered by a re-insurer, with $x$ the initial capital and $a$ a fixed constant. The ruin probability $P_a(x)$
equals
\begin{equation}
 P_a(x)=   \Prob{ \sup_{n\geq 0} \sum_{i=1}^n [B_i I(B_i < ax) - c A_i]>x},
\end{equation}
recalling that $I(\cdot)$ is the indicator function.
We assume $c\Exp{A_i} > \Exp{B_i}$ so that $P_a(x)\rightarrow 0$. We wish to understand the impact of $a$ on the convergence rate of $P_a(x)$ and aim to
 illustrate the applicability of Theorem \ref{thm-main}. To simplify our presentation,  we assume $A_i\equiv 1/c$.


\begin{proposition}
\label{cor}
If $\Prob{B_1>x} = L(x)x^{-\alpha}$ with $\alpha>1$, and $L(cx)/L(x)\rightarrow 1$ for $c>0$ as $x\rightarrow\infty$, and $1/a$ is not an integer, then there exists a constant $C>0$ such that
    \begin{equation}
    P_a(x) = (L(x)x^{-(\alpha-1)})^{\lceil 1/a \rceil} C (1+o(1)).
    \end{equation}
\end{proposition}
The intuition behind this result is that ruin is caused by $\lceil 1/a \rceil$ large claims, similar to what has been established in finite-horizon problems \cite{albrecher2017reinsurance, chen2019efficient, dombry2022hidden}.
To prove Proposition \ref{cor},
we reduce the problem to analyzing the ruin probability over a large but finite time-horizon of length $Tx$ with $T\in (0,\infty)$. To this end, we use the bounds
\begin{equation}
     P_a(x) \geq \Prob{\sup_{n \leq Tx} \sum_{i=1}^n [B_i - 1] \geq x ; \sup_{n\leq Tx} B_i < ax } =: P_{a,T} (x),
     \label{reins-lb}
\end{equation}
and
\begin{equation}
\label{reins-ub}
     P_a(x)  \leq \Prob{\sup_{n \leq Tx} \sum_{i=1}^n [B_iI(B_i<ax) - 1] \geq x} + \Prob{\sup_{n \geq Tx} \sum_{i=1}^n [B_iI(B_i<ax) - 1] \geq x}.
\end{equation}
The first term in the right-hand side of (\ref{reins-ub}) can be upper bounded by
\begin{equation}
   \Prob{\sup_{n \leq Tx} \sum_{i=1}^n [B_iI(B_i<ax) - 1] \geq x \mid \sup_{i\leq Tx} B_i < ax} =  P_{a,T}(x)/ \Prob{\sup_{i\leq Tx} B_i < ax}.
   \label{reins-ub1a}
\end{equation}
To bound the second term, define $X_i = B_i-1$, note that $B_iI(B_i<ax) - 1 \leq \min \{ X_i, ax\}$, to obtain.
\begin{equation}
    \Prob{\sup_{n \geq Tx} \sum_{i=1}^n [B_iI(B_i<ax) - 1] \geq x} \leq \Prob{\sup_{n \geq Tx} \sum_{i=1}^n \min \{ X_i, ax\} \geq x}.
\end{equation}
This expression can be bounded further by observing that
\begin{align}
& \Prob{\sup_{n \geq Tx} \sum_{i=1}^n [B_iI(B_i<ax) - 1] \geq x} \leq \Prob{\sup_{n \geq Tx} \sum_{i=1}^n \min \{ X_i, ax\} \geq x}.
 \nonumber \\
 &\leq \Prob{\sup_{n \geq Tx} \sum_{i=1}^n \min \{ X_i, ax\} \geq x, \sum_{i=1}^{Tx} \min \{ X_i, ax\} \leq Tx/2} + \Prob{\sum_{i=1}^{Tx} \min \{ X_i, ax\} \geq Tx/2} \nonumber \\
&\leq \Prob{\sup_{n \geq Tx} \sum_{i=Tx+1}^n \min \{ X_i, ax\} \geq xT/2} + \Prob{\sup_{n\geq 0} \sum_{i=1}^{n} \min \{ X_i, ax\} \geq Tx/2} \nonumber \\
  &=2 \Prob{M(ax) > Tx/2}, \label{reins-ub1b}
\end{align}
with $M(ax)$ being the protagonist of this paper, defined in (\ref{protagonist}).
Substituting (\ref{reins-ub1a}), (\ref{reins-ub1b}) into (\ref{reins-ub}), we obtain
\begin{equation}
\label{reins-ub2}
     P_a(x)  \leq P_{a,T}(x)/ \Prob{\sup_{i\leq Tx} B_i < ax} + 2\Prob{M(ax) > Tx/2}.
\end{equation}
Since $\Exp{B_i}<\infty$,
\begin{equation}
    \Prob{\sup_{i\leq Tx} B_i < ax} = \Prob{B_1<ax}^{\lfloor Tx\rfloor} \rightarrow 1.
\end{equation}
Consequently, it suffices to establish the asymptotic behavior of
$P_{a,T}(x)$ and to show that
\begin{equation}
\label{cutoff}
  \Prob{M(ax) > Tx/2} = o( P_{a,T}(x))
\end{equation}
for a suitably chosen $T$.
In Section 5.1 of \cite{rhee2016arxiv}, it is shown that, given $a$ is such that $1/a$ is not an integer, there exists a constant $C$ such that
\begin{equation}
\label{finite-time-pareto}
    P_{a,T}(x) \sim C (L(x)x^{-(\alpha-1)})^{\lceil 1/a \rceil} (1+o(1)).
\end{equation}
The constant $C$ can be expressed as
\[
C= \int_{u\in [0,T]^k, x\in (0,\infty)^k} \alpha^k \prod_i x_i^{-1-\alpha} I(\sup_{t\in [0,T]} \sum_{i=1}^k x_i I(u_i \leq t)-t \geq 1 ; \sup_{i\leq k} x_i \leq a)
du_1...du_k dx_1...dx_,
\]
with $k=\lceil 1/a\rceil$.
From this representation, it follows that $C \in (0,\infty)$ if $1/a$ is not an integer. In addition, $C$ is constant in $T$ for $T>1+a$.
Next, using (\ref{eq-main-beta}), we see that, for every $\beta \in (1,\alpha)$, and $x$ large enough,
\begin{equation}
      \Prob{M(ax) > Tx/2} \leq (ax) ^{-(\beta-1) T/(2a)},
\end{equation}
which, using (\ref{finite-time-pareto}), is $o( P_{a,T}(x))$ for $T> \frac{\alpha-1}{\beta-1}a\lceil 1/a \rceil$, concluding the proof of Proposition \ref{cor}. \\

\newpage

\noindent
{\bf Remarks.}
\begin{enumerate}
\item
To derive tail asymptotics for the ruin probability if $\Prob{X_1>x}$ is of Weibull-type, we expect that it is
possible to apply the sample path large-deviations results in \cite{bazhba2020sample} to derive tail asymptotics for $P_{a,T}(x)$, and combine them with the
bound (\ref{reins-ub2}) and Theorem \ref{thm-main2}; a fully worked-out argument would require a careful investigation of a quasi-variational problem, extending the analysis in Section 4 of \cite{bazhba2020sample}, which is beyond the scope of this study.
\item In this section, we took $A_i$ constant. The case of non-deterministic $A_i$ may be dealt with by using bounds of the form
$$\sup_{n\geq 0} \sum_{i=1}^n [B_i I(B_i < ax) - c A_i] \leq \sup_{n\geq 0} \sum_{i=1}^n [B_i I(B_i < ax) - c']+ \sup_{n\geq 0} \sum_{i=1}^n [c'- cA_i],$$
with $c'\in (\Exp{B_i}, c\Exp{A_i})$. The second term on the right-hand side of this inequality is independent of the first term. It has a moment generation function that is finite in a neighborhood of the origin
and therefore does not impact the asymptotic behavior of $P_a(x)$.
\end{enumerate}


\bibliographystyle{plain}


\appendix
\section{Bounding Taylor approximation errors}

In the proof of Theorem \ref{thm-main}, the following quantities appear:
\begin{equation}
    ME_\delta = \max_{s\geq 0} \frac{1-(1+s)e^{-s}}{s^\delta}, \hspace{1cm} \delta \in [0,2],
    \label{a1}
\end{equation}
\begin{equation}
    MG_\delta = \max_{u\geq 0} \frac{e^{-u}-1+u}{u^\delta}, \hspace{1cm} \delta \in [1,2].
    \label{a2}
\end{equation}
We analyze these quantities using the machinery developed in \cite{Janssen2024}, which contains a framework to bound Taylor polynomial approximation errors for the exponential function in the presence of a power weight function. For the case that the degree $n$ of the Taylor polynomial equals 1, this leads to consideration of quantities $ME_{1,\delta}$ and $MG_{1,\delta}$ that agree with (\ref{a1}) and (\ref{a2}), from which the degree $n=1$ has been omitted for brevity. In this appendix, we summarize the results of \cite{Janssen2024} for the case $n=1$ that are relevant to this paper. We have
\begin{equation}
\label{a3}
    ME_0=1, ME_2= \frac 12, MG_1=1, MG_2 = \frac 12.
\end{equation}
Furthermore, $ME_\delta$ is a continuous, log-convex function of $\delta\in [0,2]$, and $MG_\delta$ is a continuous, log-convex function of $\delta \in [1,2]$. For $\delta \in (0,2)$, the maximum $ME_\delta$ of
\begin{equation}
\label{a4}
    E_\delta(s) = \frac{1-(1+s)e^{-s}}{s^\delta}
\end{equation}
over $s\geq 0$ is assumed by the unique positive solution $s=s(\delta)$ of the equation
\begin{equation}
\label{a5}
    e^s = 1+s+\frac 1\delta s^2.
\end{equation}
For $\delta\in (0,1)$, the maximum $MG_\delta$ of
\begin{equation}
\label{a6}
G_\delta(u) = \frac{e^{-u}-1+u}{u^\delta}
\end{equation}
over $u\geq 0$ is assumed by the unique positive solution $u=u(\delta)$ of the equation
\begin{equation}
\label{a7}
    e^{-u} = 1 - \frac{\delta u}{u+\delta}.
\end{equation}
An upper bound for $ME_\delta, \delta \in [0,2]$ follows from Proposition 2 in \cite{Janssen2024} in the form
\begin{equation}
\label{a8}
    F_2(2-\delta); \hspace{1cm} F_2(s) = \frac 1\delta \frac{s^{2-\delta}}{1+s+\frac 1\delta s^2}, \hspace{1cm} s\geq 0.
\end{equation}
Thus, we have
\begin{equation}
    \label{a9}
    ME_\delta \leq \frac{(2-\delta)^{2-\delta}}{4-\delta}=:UE_\delta, \hspace{1cm} \delta\in [0,2].
\end{equation}
An upper bound for $MG_2, \delta \in [1,2]$ is found using the approach given at the end of Subsection 2.2 of \cite{Janssen2024} in the form
\begin{equation}
    \label{a10}
    H\left(\frac{(2-\delta)\delta}{\delta-1}\right); \hspace{1cm} H(u) = \frac{u^{2-\delta}}{u+\delta}, \hspace{1cm} u\geq 0.
\end{equation}
Thus, we have
\begin{equation}
\label{a11}
    MG_\delta \leq \left(\frac{\delta-1}{\delta}\right)^{\delta-1}(2-\delta)^{2-\delta}=: UG_\delta, \hspace{1cm} \delta \in [1,2].
\end{equation}
The bounds in (\ref{a9}) and (\ref{a11}) are relatively sharp, see Figure 1 and Figure 2 below {\bf to be added}. Less sharp, but simpler (and still effective) bounds are
\begin{equation}
    \label{a12}
    ME_\delta \leq \frac 1{1+\frac 12 \delta}, \hspace{0.5cm} \delta \in [0,2]; \hspace{1cm} MG_\delta \leq \frac 1\delta, \hspace{0.5cm} \delta \in [1,2].
\end{equation}
The bound for $MG_\delta$ also appears in the proof of Lemma 1.4 of \cite{Nagaev1979}. We show for completeness that the bounds in (\ref{a12}) follow from the two bounds in (\ref{a9}) and (\ref{a11}), respectively.

The first bound in (\ref{a12}) follows via (\ref{a9}) from the inequality
\begin{equation}
    \label{a13}
    f_E(\delta) := \ln (1+\frac{1}{2}\delta) + (2-\delta) \ln (2-\delta) - \ln (4-\delta) \leq 0, \hspace{0.5cm} \delta \in [0,2].
\end{equation}
To establish (\ref{a13}), we observe that the function $f_E(\delta)$ is continuous in $\delta \in [0,2]$ and that $f_E(0)=0=f_E(2)$. Furthermore, a simple computation yields 
\begin{equation}
\label{a14}
    f_E''(\delta) = \frac {-1}{(2+\delta)^2}+ \frac 1{2-\delta}+ \frac 1{(4-\delta)^2} > -\frac 14+ \frac 12 + \frac 1{16}>0, \delta \in (0,2).
\end{equation}
Hence, $f_E(\delta)$ is convex in $\delta \in (0,2)$, and so we get (\ref{a13}) from $f_E(0)=0=f_E(2)$ and continuity of $f_E(\delta), \delta \in [0,2]$.

The second bound in (\ref{a12}) follows via (\ref{a11}) from the inequality
\begin{equation}
\label{a15}
f_G(\delta):= (2-\delta) \ln \delta + (\delta-1) \ln (\delta-1) + (2-\delta) \ln (2-\delta)\leq 0, \hspace{0.5cm} \delta \in [1,2].
\end{equation}
To establish (\ref{a15}), we observe that the function $f_G(\delta)$ is continuous in $\delta \in [1,2]$ and that $f_G(1)=0=f_G(2)$. Furthermore, a simple computation yields
\begin{equation}
\label{a16}
    f_G''(\delta) = -\frac 2{\delta^2}-\frac 1\delta + \frac 1{\delta-1}+\frac 1{2-\delta}, \hspace{0.5cm} \delta\in (1,2).
\end{equation}
Consequently,
\begin{equation}
    \label{a17}
    \delta^2(\delta-1)(2-\delta) f_G''(\delta) = -(\delta-1)(2-\delta)(2+\delta)+\delta^2> -\frac 14 4 + 1^2 = 0, \hspace{0.5cm} \delta \in (1,2).
\end{equation}
Hence, $f_G(\delta)$ is convex in $\delta \in (1,2)$, and therefore we get (\ref{a15}) from $f_G(1)=0=f_G(2)$ and continuity of $f_G(\delta), \delta \in [1,2]$.

In Figure 1, we show a plot of $ME_\delta$, together with plots of the upper bounds $UE_\delta$ in (\ref{a9}) and $(1+\frac 12 \delta)^{-1}$ in (\ref{a12}), as a function of $\delta \in [0,2]$. The value of $ME_\delta$ is obtained as $E_\delta (s(\delta))$, where $E_\delta$ is given in (\ref{a4}) and the unique solution $s=s(\delta)$ of the equation in (\ref{a5}) is computed using Newton's method (starting value: $s^{(0)} = \ln (2/\delta)+(2-\delta)$, see \cite{Janssen2024}, (37)).

In Figure 2, we show a plot of $MG_\delta$, together with plots of the upper bounds $UG_\delta$ in (\ref{a11}) and $\delta^{-1}$ in (\ref{a12}), as a function of $\delta \in [1,2]$. The value of $MG_\delta$ is obtained as $G_\delta(u(\delta))$, where $G_\delta$ is given in (\ref{a6}) and the unique solution $u=u(\delta)$ of the equation in (\ref{a7}) is computed using Newton's method (starting value: $u^{(0)}=(2-\delta ) \delta/ (\delta-1)- \ln (\delta-1)$, see \cite{Janssen2024}, (275)).

\begin{figure} 
    \centering
    \includegraphics[width=0.7\linewidth]{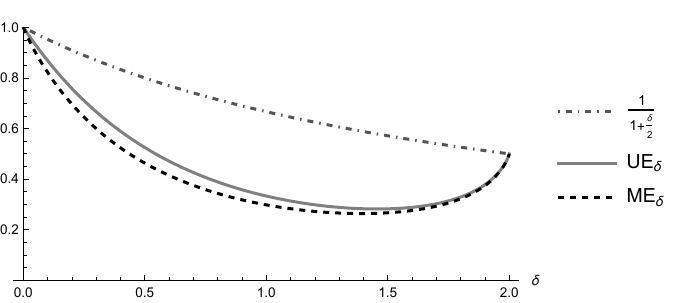}
    \caption{$ME_\delta$, $UE_\delta$, $1/(1+\delta/2)$ as function of $\delta \in (0,2)$.}
\end{figure}

\begin{figure} 
    \centering
    \includegraphics[width=0.7\linewidth]{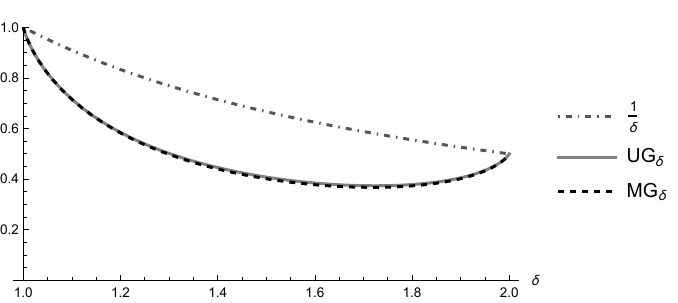}
    \caption{$MG_\delta$, $UG_\delta$, $1/\delta$ as function of $\delta \in (1,2)$.}
\end{figure}

\end{document}